\documentclass[a4paper]{article}

\usepackage{inputenc}
\usepackage[english]{babel}
\usepackage{graphicx, here, subfigure}
\usepackage{amsmath}
\usepackage{color}
\usepackage{amssymb}
\usepackage{cite}
\usepackage{mdwmath}
\usepackage{multirow}
\usepackage{mdwtab}
\usepackage{verbatim}
\usepackage{enumerate}
\usepackage{amsfonts}
\usepackage{eurosym}

\usepackage{amsthm}
\usepackage{fancyheadings}
\usepackage{fonte}

\usepackage{tikz}
\usetikzlibrary{positioning,chains,scopes,arrows,backgrounds,fit,mindmap,trees}

\newcommand{\RRR}{{\mathrm I\! \textsc{R} }} 
\newcommand{\K}{{\kappa }} 
\newcommand{\KKK}{{\mathcal{K}}} 
\newcommand{\VVV}{{\mathcal{V}}} 
\newcommand{\TTT}{{\mathcal{T}}} 
\newcommand{\XX}{{X^\# }} 
\newcommand{\PP}{{P^\# }}

\newtheorem{proposition}{\textsc{ \textbf{Proposition} }}[section]

\definecolor{bluerenault}{rgb}{0,0.555,0.880}
\definecolor{greenrenault}{RGB}{100,133,33}
\definecolor{yellowrenault}{rgb}{1,0.853,.255}
\definecolor{redrenault}{RGB}{240,41,37}

\title{Reachability of Delayed Hybrid Systems Using Level-set Methods}

\author{G. Granato\thanks{Renault SAS, Advanced Electronics Division, 
        TCR RUC T 65, 78286 Guyancourt Cedex, France
        ({\tt giovanni.granato@renault.com}).}}

\begin{document}
	
\thispagestyle{fancy}
\lhead{{\ninerm \bf  EngOpt 2012 - International Conference on Engineering Optimization} \\
{\eightrm Rio de Janeiro, Brazil, 1-5 July 2012.}}
\setlength{\headrulewidth}{0pt}
\begin{center}
{\large {\elevenrm \bf Reachability of Delayed Hybrid Systems Using Level-set Methods}}\\~\\
{\tenrm \bf {G. Granato }}
\\~\\
{\ninerm
Renault SAS, Guyancourt, France. \\
\'Ecole Nationale Sup\'erieure de Techniques Avanc\'ees,
        Unit\'e de Math\'ematiques Appliqu\'ees, Paris, France.\\
        giovanni.granato@renault.com.}
\end{center}
	
	\begin{abstract}
		This study proposes an algorithm to synthesize controllers for the power management on board hybrid vehicles that allows the vehicle to reach its maximum range along a given route.
		The algorithm stems from a level-set approach that computes the reachable set of the system, i.e., the collection of states reachable from a certain initial condition via the computation of the value function of an optimal control problem.
		The discrete-time vehicle model is one of a particular class of hybrid vehicles, namely, range extender electric vehicles (REEV).
		This kind of hybridization departures from a full electric vehicle that has an additional module -- the range extender (RE) -- as an extra energy source in addition to its main energy source -- a high voltage battery.
		As an important feature, our model allows for the switching on and off of the range extender and includes a decision lag constraint, i.e., imposes two consecutive switches to be separated by a positive time interval.
		The approach consists in the introduction of an adequate optimal control problem with lag constraints on the switch control whose value function allows a characterization of the reachable set.
		The value function is in turn characterized by a dynamic programming algorithm.
		This algorithm is implemented and some numerical examples are presented.
	\end{abstract}

	\pagestyle{myheadings}
	\thispagestyle{plain}
	\markboth{G. GRANATO}{REACHABILITY OF DELAYED HYBRID SYSTEMS USING LEVEL-SET METHODS}

	\section{Introduction}

	Electrified automobile power chain technology faces a rapid development.	
	The interest for this technology arises as an increasing number of automakers wish to adopt vehicles with an electrification of the power chain as a viable solution for reducing greenhouse gas emissions worldwide to meet stringer regulative legislation and consumers' demand.
	There is an effort of major constructors to deploy fully electric vehicles (EVs) as early as possible in some car market segments. 
	This work focuses on a specific class of hybrid electric vehicles, namely, range extender electric vehicles (REEV).
	The range extender electric vehicle model consists of a small dimensioned internal combustion engine (ICE) as an extra energy source acting as a range extender in addition to the main energy source, namely, a DC battery. 
	The two energy sources are assembled as in a series hybrid electric vehicle, which means that the ICE is not connected mechanically to the drive shaft. 
	Instead, a generator transforms the mechanical energy produced in the ICE into electric energy that can be directed towards the electric motor or the battery. 
	Because of the relative low power of the ICE, this configuration does not allow the vehicle to rely solely on the range extender and thus the vehicle cannot assure traction if there is no energy in the battery. 
	Although our controllers have to take into account the particularities of this class of hybrid vehicles, the formulation developed in this study is general enough and may be adapted to other architectures, including series and parallel hybrid electric vehicles.

	The goal of this work is compute a way of operating the range extender to enable the vehicle to reach its maximum range.
	The REEV is modeled as a discrete-time dynamical system in which the state vector represents the energy capacities of the two different energy sources.	
	In this setting, the study aims at finding the control sequence of the two energy sources that allows the vehicle to reach the furthest possible point of a given route or to reach the driver selected energy state.

	To the extent of our knowledge, there are no studies considering the driving range of the vehicle as an optimization criterion.
	Usually, in the context of power management strategies the optimality of a strategy is defined in terms of the total fuel consumption during a driving cycle.
	We refer here to \cite{sciarretta2004optimal,kirschbaum2002determination,caux2010line,brahma2000optimal,lin2003power}.
	Different techniques are used in each of these works, as well as different assumption on whether one has access to \emph{a priori} driving information.
	A common feature of the works above is the fact that the power management synthesis uses a predetermined \emph{driving range}, typically given as one considers a particular driving cycle.

	The main idea in our work is to have no \emph{a priori} information about the vehicle driving range.
	A typical setting is when the driving cycle is too long to be completed with the current energy on-board the vehicle.
	In the case where no power management strategy can ensure the cycle completion, the utilization of a fuel consumption criteria does not seem to lead to well-posed posed optimal control problem.
	For instance, dynamic programming-based algorithms cannot be implemented, since they rely on a backward computation of the value function starting from the "last" reachable point on the route which may be unknown.	
	This work proposes a method to synthesize a power management strategy that maximizes the driving range on a cycle supposedly too long to be completed.

	This paper is organized as follows: firstly, it describes the REEV model motivating the study and states the associated hybrid optimal control problem.
	Then, the reachable set and the value function are defined and a dynamic programming principle for the value function is obtained.
	It follows that the dynamic programming principle is used to derive an dynamic programming algorithm to compute the value function thus characterizing the reachable set of the system.
	Lastly results of numerical simulations evaluating the autonomy both of a REEV toy and a realistic model are presented.

	\section{Problem Settings}

	\subsection{Vehicle Model}
	
	This section discusses the model of the power management system of a REEV that is to be optimized and introduces the notation used throughout the document.
	A REEV is a vehicle that combines a primary power source - a HV battery - and a small dimensioned (powerwise) secondary power source - in our case an ICE. 
	The traction (or propulsion) of the vehicle is performed by an electric motor connected to the vehicle's wheels through a reduction gear. 
	Both energy sources can supply some share of the necessary power to meet the driver demand. 
	Additionally, since the model considers a range-extender electric vehicle type, it cannot rely solely on the RE power to drive the vehicle.
	Thus, the model considers that the vehicle's traction capability is conditioned to the existence of some electric energy in the battery. 

	The architecture is that of a series hybrid electric vehicle, which means that the ICE is not mechanically connected  to the transmission.
	Instead, a generator transforms the mechanical energy produced in the ICE into electric current that can be directed towards the electric motor or charge the battery. 

	The ICE state is controlled by a discrete sequence of switching orders.
	An important feature of the REEV model is a \emph{time lag} $\delta > 0$ imposed between two consecutive decisions times. 		
	From the physical viewpoint, this assumption incorporates the fact that frequent switching of the RE is undesirable in order to avoid mechanical wear off and acoustic nuisance for the driver. 
	This lag has an order of magnitude of about $120$s. 

	The controlled variables are the state of charge (SOC) of the battery, the remaining fuel on the ICE reservoir and the ICE state on or off.

	This setting considers that the optimal power management strategy is synthesized for a specific trip.
	In this direction, we assume that information about the vehicle's future power demand is known.
	In particular, it is assumed that the vehicle future speed, i.e. the speed profile of the vehicle, as well as the relief of the route followed by the vehicle are known.
	As a consequence, the problem can be stated in a deterministic setting.
	This may seem a very strong assumption.
	However, range extender electric vehicles are expected to be equipped with embedded navigation systems (NAV) that actually have information about expected speeds and topological data.
	This information can be exploited in order to construct an expected speed profile with a certain degree of confidence.
	
	In the following of the paper, index $k=1,\cdots,\K$ shall index a position in the route.
	Since the vehicle speed is assumed to be known, on can easily obtain, from a known location $k$, the time to attain it.
	More precisely, denote by $d_k$ the distance between points $k-1$ and $k$ and $\sigma_k>0$ the (constant) vehicle speed between points $k-1$ and $k$.
	The total time to drive between $k-1$ and $k$ is then simply $d_k/\sigma_k$.
	
	The state vector is composed of the vehicle energy state and the ICE state.
	The vehicle energy state is a two-dimensional vector  $y \in X = \RRR^2$, where $y=(y_1,y_2)$ denotes the state of charge of the battery and the fuel available in the range extender module respectively.
	Each of these quantities are the image of the remaining energy in the battery and in the RE respectively.
Moreover, because energy capacities of the battery and the fuel tank are limited (and normalized between $1$ -- full reservoir -- and $0$ -- empty reservoir), the set of admissible states is defined as $K = [0,1]^2$.
	Remark that we define the state space as all the plane $X = \RRR^2$ and introduce a set of admissible states $K$.
	The set $K$ represents the "physically admissible" SOC and fuel levels, considering that they are normalized capacities.
	This is a modeling choice that later enables us to deal easily with the state constraints in the dynamic programming principle (in opposition to imposing "hard" constraints on the state space like $X=[0,1]^2$).
	The RE state is denoted by $q \in Q = \{0,1\}$ and indicates whether the RE is off ($q=0$) or on ($q=1$).

	The control variables are the output power of the ICE and the order to switch the ICE's state (on or off).
	The first is referred to as the continuous control and the latter as the discrete (or switch) control.
	The continuous control is applied at all points $k=0,\cdots,\K-1$ of the route and is denoted $u_k \in U(q_k)$.
	The idea is that at each point $k$, also called the decision stages, the controller is able to control the power split of the vehicle.
	Indeed, in a series architecture the power split is controlled only by setting the output power of the range extender.
	That is because, for a given required power at the wheel $P^*$, imposed by the driver, once the range extender output power $u$ is chosen the remaining required power $P^*-u$ is automatically supplied by the battery.

	The discrete control $w$ is a sequence of switching decisions $w =\{ (w_1,s_1), \cdots, (w_N,s_N)\}$, where  for all $j = 1,\cdots,N$, $s_i \in \{1,\cdots,\K-1 \}$ and $w_j \in W(q_{s_j}) \subset Q$.
	The sequence of switch controls $\{w_j\}_{j=1}^N$ (designating the new mode of operation of the ICE) is associated with the sequence of switching positions $\{s_j\}_{j=1}^N$ so that each decision $w_j$ is exerted at node $s_j$.
	The set $W(q)$ is a set value function of the actual discrete state and contains all possible new discrete states.
	In the vehicle application, where there are only two discrete modes of operation (on or off), this set is simply $W(1) = \{0\}$ and $W(0) = \{1\}$.
	The switching decision is introduced so one is able to stop the range extender whenever there are large idle periods, thus saving fuel.
	More precisely, if the ICE is running and no output power is required, there is still an idle fuel consumption, as energy is required to keep the engine running to avoid stalling.
	In the case of extended idle periods, it may be worth to stop the ICE and pay a fuel overhead to restart the engine whenever needed.

	The switch lag condition imposes that two consecutive switch decisions must be separated by a time interval of $\delta >0$, i.e.
	\begin{equation*}
		\sum_{i = s_j+1}^{s_{j+1}} \frac{d_i}{\sigma_i}\geq  \delta.
	\end{equation*}
	
	Given a discrete control sequence $w$ and a discrete state $q_j$, the continuous control steers the state according to
	\begin{eqnarray*}
		y_{k-1} = f_k(y_k,u_{k-1},q_j) 
	\end{eqnarray*}	
	and at isolated points $\{s_j\}$, the discrete state changes according to
	\begin{equation*}
		q_{j-1} = g(q_j,w_j).
	\end{equation*}
	These expressions mean that: given a discrete control sequence $w$, let $q_j$ be a mode of operation activated at some point $s_{j+1}$.
	Then, the dynamic $f_\cdot(\cdot,\cdot,q_j)$ operates on the state between points $k=s_{j+1},\cdots,s_j$ and the trajectory of the discrete state is simply $q_k=q_j$ for $k$ between $s_{j+1}$ and $s_j$.
	When point $s_j$ is reached, the discrete dynamics undergoes a switch associated with the switch control $w_j$ and changes to $q_{j-1} = g(q_j,w_j)$, activating the new mode of operation given by $f_\cdot(\cdot,\cdot,q_{j-1})$ .
	Observe moreover that the process is modeled and described in backward fashion.

	At this point, we introduce the notion of hybrid control and admissible hybrid control.
	The hybrid control is introduced for notational convenience.
	It is a control regrouping both the range extender power $u$ and the switch decision $w$.
	More in detail, given a discrete control sequence $w =\{ (w_1,s_1), \cdots, (w_N,s_N)\}$ and a continuous control sequence $u = (u_1,\cdots,u_{\K-1})$, the hybrid control is denoted by $a= (u,w)$ and groups the continuous and discrete controls.
	Denote by $A$ the set of all hybrid controls.
	As mentioned,  the hybrid control $a$ must  contain controls $(u,w)$ such that $u_k \in U(q_k)$ and $w_j \in W(q_{s_j}) \text{ for } j=1,\cdots,N$.
	However, the admissibility of each discrete and continuous control does not imply the admissibility of the hybrid control without additional structural conditions.
	By structural conditions we refer to the fact that the decision nodes $s_j$ must be increasing with $j$ and that the associated decision times must respect the lag condition.
	The structure imposed on the admissible hybrid controls $a\in A^{\text{adm}}$ allows the system to engender admissible trajectories.
	An admissible trajectory can be regarded as a "physically" admissible trajectory.
	Such trajectories correspond to the evolution of the SOC and fuel levels whenever the system is steered by an admissible control sequence.	
	However, the hybrid control admissibility condition is not well adapted to a dynamic programming principle formulation, needed later on.
	In order to include the admissibility condition in the optimal control problem in a more suitable form, we introduce a new state variable $\pi$.
	Recall that the decision lag condition implies that new switch orders are not available until a time $\delta$ after the latest switch.
	The new variable is constructed such that, at a given position $k$, the value $\pi_k$ measures the time since the last switch.
	The idea is to translate the structural conditions into conditions satisfied by this new variable and treat them as state constraints.
	Thus, if $\pi_k<\delta$ all switch decisions are blocked and if, conversely, $\pi_k \geq \delta$ the system is free to switch.
	For that reason, this variable can be seen as a switch lock.
	Now, given $\K >0$, $k = 1,\cdots,\K$, a discrete control sequence $w$, final conditions $x \in X$ and $q\in Q$, the hybrid dynamical system is written as 
	\begin{eqnarray}
		y_{k-1} 	&=& f_k(y_k,u_{k-1},q_j), ~~ k=\K,\cdots,1, ~~~~~~~~y_ \K 	= x 	\label{hybDynamicA} \\
		q_{j-1} 	&=& g(q_j,w_j), ~~ j=1,\cdots,N, ~~~~~~~~~~~~~~~~q_N = q.	\label{hybDynamicB} \\
		\pi^w_k &=& \pi_k = \left\{ \begin{array}{lll}
						\delta + \sum_{i=1}^k \frac{d_i}{\sigma_i} 	& \text{if} & k < s_1 \\
						\inf \limits_{\substack{s_j \leq k} }\sum_{i=s_j}^k \frac{d_i}{\sigma_i} - \frac{d_{s_j}}{\sigma_{s_j}} 	& \text{if} & k \geq s_1 \\
					\end{array} \right.	\label{pDynamic}
	\end{eqnarray}

	Indeed, once the discrete control is given, the trajectory $\pi_k$ can be determined.
	Proceeding with the idea of adapting the admissibility condition we wish to consider $\pi_\K=p$, with $p \in P := (0,T]$, where $T$ is the total travel time, the final value of the switch lock variable trajectory and impose the delay condition under the form $\pi_{s_j-1} \geq \delta$ for all $s_j$ (notice that $\pi_{s_j}=0$ by construction).
	Proceeding as such, these conditions suffice to define an admissible discrete control set.
	So, while optimizing with respect to admissible functions, one needs only look within the set of hybrid controls that engender a trajectory $\pi_k$ with an appropriate structure.
	Denote the solutions of \eqref{hybDynamicA}-\eqref{pDynamic} with final conditions $x,q,p$ by $y_{x,q,p;\K}$,$q_{x,q,p;\K}$ and $\pi_{x,q,p;\K}$.
	Given $\K >0$, $x \in X$, $q \in Q$ and $p \in P$, define the admissible trajectory set $ S_{(0,\K)}^{x,q,p} $ as
	\begin{eqnarray}	
		\nonumber S_{(0,\K)}^{x,q,p} &=& \{y_{x,q,p;\K} ~|~ a= (u,\{w_j,s_j \}_{j =1}^{N})  \in A, ~y_{x,q,p;\K} \text{ solution of \eqref{hybDynamicA}-\eqref{hybDynamicB}}, ~~~~~~~~~~ \\
					&&	~~~~~~~~\pi_{x,q,p;\K} \text{ solution of \eqref{pDynamic}}, ~\pi_\K = p, ~\pi_{s_j-1} \geq \delta, ~ j=1,\cdots,N \}.
		\label{admTrajSet}
	\end{eqnarray}

	\subsection{Reachability Problem and Autonomy}

	Let $X_0 \subset X$ be the set of allowed initial states, i.e. the set of states from which the system \eqref{hybDynamicA}-\eqref{hybDynamicB} is allowed to start.
	Define the reachable set as the set of all points attainable by $y$ after driving until point $k$ starting within the set of allowed initial states $X_0$ to be
	\begin{eqnarray}
	\nonumber	R^{X^0}_k 	&=&\left \{ x \mid \exists (q,p) \in Q \times P,~   y_{x,q,p;k} \in S^{x,q,p}_{(0,k)}, \right. 			\\
					&& ~~~~~~~~~~~~~~~~~~\left. (y_{x,q,p;k})_0 \in X^0, \text{ and } (y_{x,q,p;k})_\theta \in K,~ \theta = 0,\cdots,k \right\}
	\label{reachForward}
	\end{eqnarray}

	In other words, the reachable set $R^{X^0}_k$ contains the values of the final energy state $(y_{x,q;k})_k$, regardless of $q,p$ at point $k$, for all admissible trajectories -- i.e., trajectories obtained through an admissible hybrid control -- starting within the set of possible initial states $X_0$ that never leave set $K$. 
	In particular, the information contained in \eqref{reachForward} allows one to determine the first point where the reachable set is empty.
	More precisely, given $X_0 \subset X$, define $k^* \geq 0$ to be
	\begin{equation}
		k^* = \inf \{ k \geq 0 \mid R^{X^0}_k \subset \emptyset \}.
	\label{hybridAutonomy}
	\end{equation}
	The position \eqref{hybridAutonomy} is identified as the \emph{autonomy} of the hybrid system \eqref{hybDynamicA}-\eqref{hybDynamicB}.
	Indeed, one can readily see that if no more admissible energy states are attainable after $k^*$, any admissible trajectory must leave set $K$ of admissible values beyond this time.
	Therefore, $k^*$ is interpreted as the last point in the path in which there is some usable energy on-board the vehicle.

	In the case of range extender electric vehicles, the vehicle stops operating in nominal mode when the battery is depleted, even if there is some fuel left in the tank.
	If the driving conditions allows all the fuel to be consumed, then the last reachable energy state of the vehicle is the state $(0,0)$, where there is no charge in the battery and no fuel in the tank.
	However, this is by all means the sole last reachable energetic state in all situations.
	Indeed, because the range extender is power-wise small dimensioned, it may be the case that the engine cannot consume energy quickly enough to empty the fuel tank.
	For instance, if the vehicle is driving through a steep slope, the battery may be depleted before all the fuel is consumed resulting thus in a final state $(0,y_2)$, $y_2>0$.
	The reachable set $R^{X^0}_k $ is associated with a \emph{minimum time function} $\TTT : X \to \RRR^+ \cup \{ \infty \}$, defined as

	\begin{equation} \label{minimumTimeDef}
		\TTT(x) 	= \inf \{ k \geq 0 \mid x \in R^{X^0}_k \}
	\end{equation}
	The minimum time function associated to a point $x$ in the state space is the minimum time that any admissible trajectory needs to reach $x$ departing from set $X_0$ while respecting the state constraints $K$.

	\subsection{Optimal Control Problem and Dynamic Programming Principle}

	In order to characterize the reachable set $R^{X^0}_k$ this paper follows the classic level-set approach \cite{OsherSethian88}.
	The idea is to describe \eqref{reachForward} as the negative region of a function $v$. 
	It is well known that the function $v$ can be defined as the value function of some optimal control problem. 
	In the case of system \eqref{hybDynamicA}-\eqref{hybDynamicB}, $v$ happens to be the value function of a (discrete time) hybrid optimal control problem.
	Consider a Lipschitz continuous function $\phi : X \rightarrow \RRR$ such that
	\begin{equation} \label{distance}
		\phi(x) \leq 0 \Leftrightarrow x \in X^0.
	\end{equation}
	Such a function always exists -- for instance, the signed distance function $d_{X_0}$ from the set $X^0$.

	For a given point $k\geq 0$ and hybrid state vector $(x,q,p) \in X \times Q \times P$, define the value function to be
	\begin{equation}
		v^0_k(x,q,p) = \inf_{ S_{(0,k)}^{x,q,p}} \left\{ \phi((y_{x,q,p;k})_0) ~|~ (y_{x,q,p;k})_\theta \in K,~ \forall \theta \in \{0,\cdots,k\} \right\}.
		\label{valueConst}
	\end{equation}
	Observe that \eqref{reachForward} works as a level-set to the negative part of \eqref{valueConst}.
	Indeed, since \eqref{valueConst} contains only admissible trajectories that remain in $K$, by \eqref{distance} implies that $v^0_k(x,q,p) $ is negative if and only if $(y_{x,q,p;k})_0$ is inside $X_0$, which in turn implies that $x \in R^{X^0}_k$.
	
	Now we proceed to describe the state constraints in the similar manner.
	The set $X \backslash K$ of state constraints is treated as an \emph{obstacle} in the state space, i.e., the state trajectory cannot pass through the obstacle set $K$.
	We recall that $K=[0,1]^2$ is the normalized SOC and fuel quantities.
	More precisely, define a Lipschitz continuous function $\varphi: X \rightarrow \RRR$ to be
	\begin{equation} \label{penal}
		\varphi(x) \leq 0 \Leftrightarrow x \in K.
	\end{equation}
	Then, for a given $k\geq 0$ and $(x,q,p) \in X \times Q \times P$, define a total penalization function to be
	\begin{equation*}
		J_k(x,q,p ; y) = \left( \phi((y_{x,q,p;k})_0) \bigvee \max_{\theta \in \{0,\cdots,k\}} \varphi((y_{x,q,p;k})_\theta) \right)
	\end{equation*}
	and then, the optimal value :
	\begin{equation}
		v_k(x,q,p) = \inf_{y\in S_{(0,k)}^{x,q,p}} J_k(x,q,p ; y).
		\label{valuePenal}
	\end{equation}
	Here and in the rest of the paper, we denote $a \vee b := \max (a,b)$.

	The next proposition certifies that \eqref{reachForward} is indeed a level-set of \eqref{valueConst} and \eqref{valuePenal}.	

	\begin{proposition}
	 \label{levelSet}
		Let $X_0 \subset X$.
		Define functions $\phi$ et $\varphi$ respectively by \eqref{penal} and \eqref{distance}.
		Define value functions $v^0$ et $v$ respectively by \eqref{valuePenal} and \eqref{valueConst}.
		Then, for $k>0$, the reachable set \eqref{reachForward} is given by
		\begin{equation}
			R_k^{X_0} = \left\{ x \mid \exists (q,p) \in Q \times P, ~ v_k(x,q,p) \leq 0 \right\} = \left\{ x \mid \exists (q,p) \in Q \times P, ~ v_k^0(x,q,p) \leq 0 \right\}.
		\label{levelSetEquation}
		\end{equation}

	\end{proposition}
	
	Proposition \eqref{levelSet} sets the equivalence between \eqref{reachForward} and the negative regions of \eqref{valueConst} and \eqref{valuePenal}.
	In particular, it states that it suffices to computes $v$ or $v_0$ in order to obtain information about $R^{X_0}$.
	In this sense, this paper focuses on \eqref{valuePenal}, which is associated with an optimal control problem with no state constraints.
	The value function can also be used to characterize the minimum time function. 
	In addition to the minimum function, the next proposition introduces the \emph{extended minimum time function}, which takes into account all state variables $x,q,p$.
	This extended function is needed later to the controller synthesis.
	\begin{proposition}
	 \label{minTimeAndValueFunction}
		The minimum time function $\TTT : X \to \RRR^+ \cup \{\infty\} $ is given by
		\begin{equation}		\label{minTimeAndValueFunctionEquation}
			\TTT(x) = \inf \left\{ k \mid \exists (q,p) \in Q \times P, ~ v_k(x,q,p) \leq 0 \right\}.
		\end{equation}
		The \emph{extended} minimum time function $\TTT' : X \times Q \times P \to \RRR^+ \cup \{\infty\} $ is given by
		\begin{equation}		\label{extMinTimeAndValueFunctionEquation}
			\TTT'(x,q,p) = \inf \left\{ k \mid v_k(x,q,p) \leq 0 \right\}.
		\end{equation}

	\end{proposition}
	Observe that equation \eqref{minTimeAndValueFunctionEquation} is an immediate consequence of the level set identities in proposition \ref{levelSet}.
	Relation \eqref{extMinTimeAndValueFunctionEquation} characterizing the extended minimum time function using the value function is derived in a similar fashion.

	Before stating a dynamic programming principle satisfied by \eqref{valuePenal}, we introduce some preliminary notations.
	Given $\K >0$, set $\KKK = \{1,\cdots,\K \}$, $\Omega = X \times Q \times P \times \KKK$ and denote its closure by $\overline \Omega$. 
	(We recall that, concerning the vehicle application, $X=\RRR^2, Q=\{0,1\}, P=(0,T],$ where $T$ is taken to be the total travel time).
	As such, \eqref{valuePenal} is defined as
	\begin{eqnarray*}
		v: \overline \Omega &\to& \RRR \\
		(x,q,p,k) &\mapsto& v_k(x,q,p)
	\end{eqnarray*}
	Define
	$$ \VVV(\overline\Omega) := \{ v ~|~ v: \overline\Omega \to \RRR, ~v \text{ bounded } \}.$$
	Now, define the non-local switch operator $M: \VVV(\overline \Omega) \rightarrow \VVV(\overline \Omega)$ to be  (adapted from the ideas in \cite{Zhang})
	\begin{equation*}
		(Mv_k)(x,q,p) = \inf_{\substack{w\in W(q) \\ p' \geq \delta}} v_k(x,g(w,q),p') \\
	\end{equation*}
	The action of this operator on the value function represents a switch that respects the delay constraint.
	Thus, they only make sense when $p=0$.
	The next proposition states the dynamic programming principle verified by \eqref{valuePenal}.

	\begin{proposition}
		The value function \eqref{valuePenal} satisfies the following dynamic programming principle:
		\begin{enumerate}
		\item[(i)] For k=0, 
		\begin{equation}
			v_0(x,q,p) = \phi(x) \bigvee \varphi(x),~~ \forall (x,q,p) \in X \times Q \times P,
		\label{initDPP}
		\end{equation}	
		\item[(ii)] For $p=0$, 
		\begin{equation}
				v_k(x,q,0) =(Mv_k)(x,q,0), ~~ \forall(x,q,k) \in X \times Q \times \KKK,
		\label{switchPDD}
		\end{equation}
		\item[(iii)] Otherwise,
		\begin{equation}
			v_k(x,q,p) = \inf_{u_{k-1}} \left\{ v_{k-1} \left( f_k(x,u_{k-1},q),q,p-\frac{d_k}{\sigma_k} \right)  \bigvee \varphi(x) \right\}.
		\label{eqPDD}
		\end{equation}
	\end{enumerate}
	\label{DynPrinciple}
	\end{proposition}

	The dynamic programming principle allows one to compute the value function and consequently characterize the reachable set $R_k^{X_0}$.
	Once the value function is obtained, one can synthesize a controller using the minimum time function $\TTT$.
	The computation of the value function and the subsequent controller synthesis are obtained in a classic fashion by a dynamic programming algorithm.
	The algorithm is the projection of the dynamic programming principle (cf. proposition \ref{DynPrinciple}) over a grid $\XX \times Q \times \PP$ of the discretized state space.
	The dynamic programming algorithm outputs a controller $(u^*,w^*)$ that controls the range extender power output and its mode of operation to reach the destination point after a time $\kappa$ such that the final state is $(x',q',p')$ departing from a state in $X_0$.
	Additionally, the trajectory remains inside set $K$, which  in our case amounts to respect the fuel and SOC capacities.
	One can observe that the extended minimum time function plays a major role in the synthesis of the controller.
	Indeed, the value function numerical values have no physical meaning since the obstacle and target functions $\phi$ and $\varphi$ are arbitrary.
	Only the sign of the value function is needed to determine the reachable set.
	Instead, the physical information is encoded in the minimum time function.

	\section{Numerical Simulations}

	Proposition \ref{DynPrinciple} provides the expressions that allow the computation through a dynamic programming algorithm.
	Once the value function obtained, one can use \eqref{levelSetEquation} to characterize the reachable set and \eqref{hybridAutonomy} to obtain the autonomy of the system.
	First, a simple vehicle model is used in the simulations.
	In addition, a constant vehicle speed is assumed during this first simulation phase.
	This model allows an analytic evaluation of the autonomy and is suitable for an a posteriori verification of the results.
	After the algorithm is tested using this simple configuration, we introduce a realistic vehicle model and a realistic speed profile.

	In this simple setting, we assume that the vehicle follows a constant speed of $\sigma \equiv 1$m/s and the distance between any two points $k-1$ and $k$ is $d_k=1$m.
	We consider that $X = \RRR^2$ and $Q=\{0,1\}$.
	The switch dynamics are given simply by $g(w,q) = w$
	ensuring that $w$ imposes the new mode of operation at every switching time.
	The continuous state dynamic models a simple REEV charge behavior.
	The evolution of the energetic state $(x,y) \in \RRR^2$ -- respectively, the battery SOC and the range extender fuel -- is given by $f(u,q) = ( -a_x + qu , -q(a_y+u))$, where $a_x,a_y >0$ are constant depletion rates of the battery electric energy and the reservoir fuel (whenever the RE is on), respectively.
	The coefficient $a_x$ models a non-negative power requirement throughout the driving cycle and therefore, a non-negative SOC consumption.
	However, the net SOC variation can be made zero or even positive (i.e., charge the battery) by the action of the output power of the range extender $u$.
	To do so, the range extender must be running and fuel must be consumed.
	The factor $a_y$ represents the fuel consumed to generate the necessary power to counter the inertia of moving parts of the range extender while $u$ is the effective share of power transfered to the electric motor.
	In order to simplify the analysis we consider the conversion factors from units of power to units of (normalized) fuel mass and percentage of battery SOC to be unitary.
	
	The engine maximum power is denoted by $u_{\text{max}}$ and the minimum power is set to be $0$.
	
	Considering this simple dynamics, an exact autonomy of the system can be evaluated analytically.
	Since there are no non-linearities in the model, an optimal strategy to reach the furthest point is to empty completely the fuel as soon as possible, thus obtaining the maximum amount of energy available.
	Given initial conditions $(x_0,y_0)$ the shortest time to empty the fuel reservoir is given by $t^* = y_0 / (a_y+u_{\text{max}})$.
	The SOC evaluated at this instant is given by $x(t^*) = x(0) -t^*(a_x - u_{\text{max}})$. 
	If $x(t^*)\leq 0$, it means the fuel cannot be consumed fast enough before the battery is depleted.  
	This condition can be expressed in terms of the parameters of the model as $x_0(a_y+u_{\text{max}})  \leq y_0 (a_x -u_{\text{max}})$.
	In this case, the autonomy is given by $T^0 = \frac{x_0}{a_x-u_{\text{max}}}$.
	Time $T^0$ is the time needed to spend all the energy on board, using the range extender at full power.
	In the other hand, if $x(t^*) > 0$, all the fuel in the range extender can be consumed and there is still some time until the battery is emptied.
	In this case, the autonomy of the system given by $T^1 = \frac{x_0 + u_{\text{max}}t^*}{a_x}$.
	\begin{figure} \label{figuresAutonomyReach}
		\includegraphics[width=1\columnwidth]{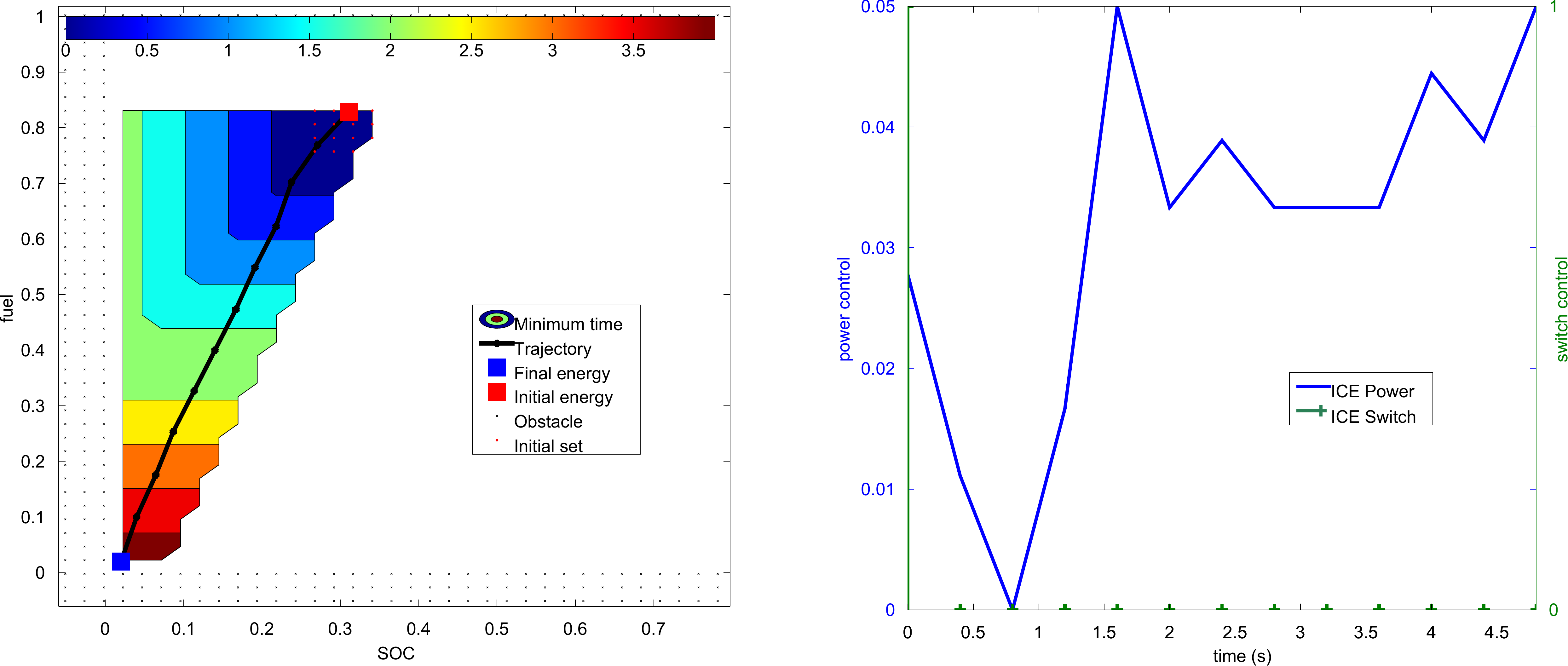} 
	\caption{Figure showing the trajectory of the energy state from a valid initial state until the last reachable point.
		The plot on the left groups the contour levels of the minimum time function, the obstacle and initial sets.
		On the right, the figure displays the range extender power control $u$ and the switch sequence $w$ controlling the ICE state.
		The synthesis is obtained using step sizes of $\Delta x = 0.0245, \Delta p = 0.4211, \Delta t = 0.4$. }
	\end{figure}

	Figure \ref{figuresAutonomyReach} illustrates a controller synthesized for a given initial condition.
	One can readily see that the reachable set, and in particular, the trajectory remain inside the authorized region of the state space $K = [0,1]^2$.
	Also, the controllers respect the delay condition $\delta=1$ between two consecutive switches.
	The initial set is taken to be a ball around $(x_0,y_0)$ in each case for numerical reasons.

	Table \ref{convSimulationsDX} presents the error between the theoretical and computed autonomy as well as the CPU time required for two initial conditions in the instance using $a_x=0.10, ~ a_y=0.15, ~u_{\text{max}}=0.07$ and $\delta = 1$.
	One can observe the convergence of the algorithm towards the theoretical value as $\Delta x$ approaches $0$.
	Empirically, a negligible sensibility with respect to $\Delta p$ is observed.

	\begin{table} \label{convSimulationsDX}
	\centering
		\begin{tabular}{|c|c|c|c|} 
		\hline
		$(x_0,y_0)$ 			&$\Delta x $ 	& $\varepsilon$ 	&  CPU\footnotemark[1]  running time(s)  \\ \hline
		\multirow{5}{*}{(0.5,0.5)}	&$0.05$ 	& $3.291$	& $2.48$		\\ 
						&$0.04$ 	& $1.154$	& $2.89$		\\ 
						&$0.03$ 	& $0.218$	& $4.21$ 		\\ 
						&$0.02$	& $0.081$	& $8.55$ 		\\ \hline
		\end{tabular}
		\caption{Convergence results and running times. Fixed $\Delta t = 0.4$ and $\Delta p = 0.5$.
			Instance considering $a_x=0.10, ~ a_y=0.15, ~u_{\text{max}}=0.07$}
	\end{table}
		\footnotetext[1]{Intel Xeon E5504 @ $2\times 2.00$GHz, $2.99$Gb RAM.}

	These remarks about the behavior of the algorithm allow us to proceed to the study of the reachable set in the realistic vehicle model.
	The procedure for implementing the simulations is as follows.
	We consider a driving profile consisting of $\kappa >0$ nodes and $\kappa-1$ links.
	At each node $k$, entry to the link $k$ there is an associated link speed $\sigma_k$ and the length of the link $d_k$.
	Once $\sigma_k,d_k$ are known, one can fix the time steps $\Delta t_k$ used in the algorithm.
	Notice that the time steps depend on the link, the problem being non-autonomous in this case.

\begin{figure} \label{vehicleAutonomyReach}
	\begin{center}
		\includegraphics[width=0.8\columnwidth]{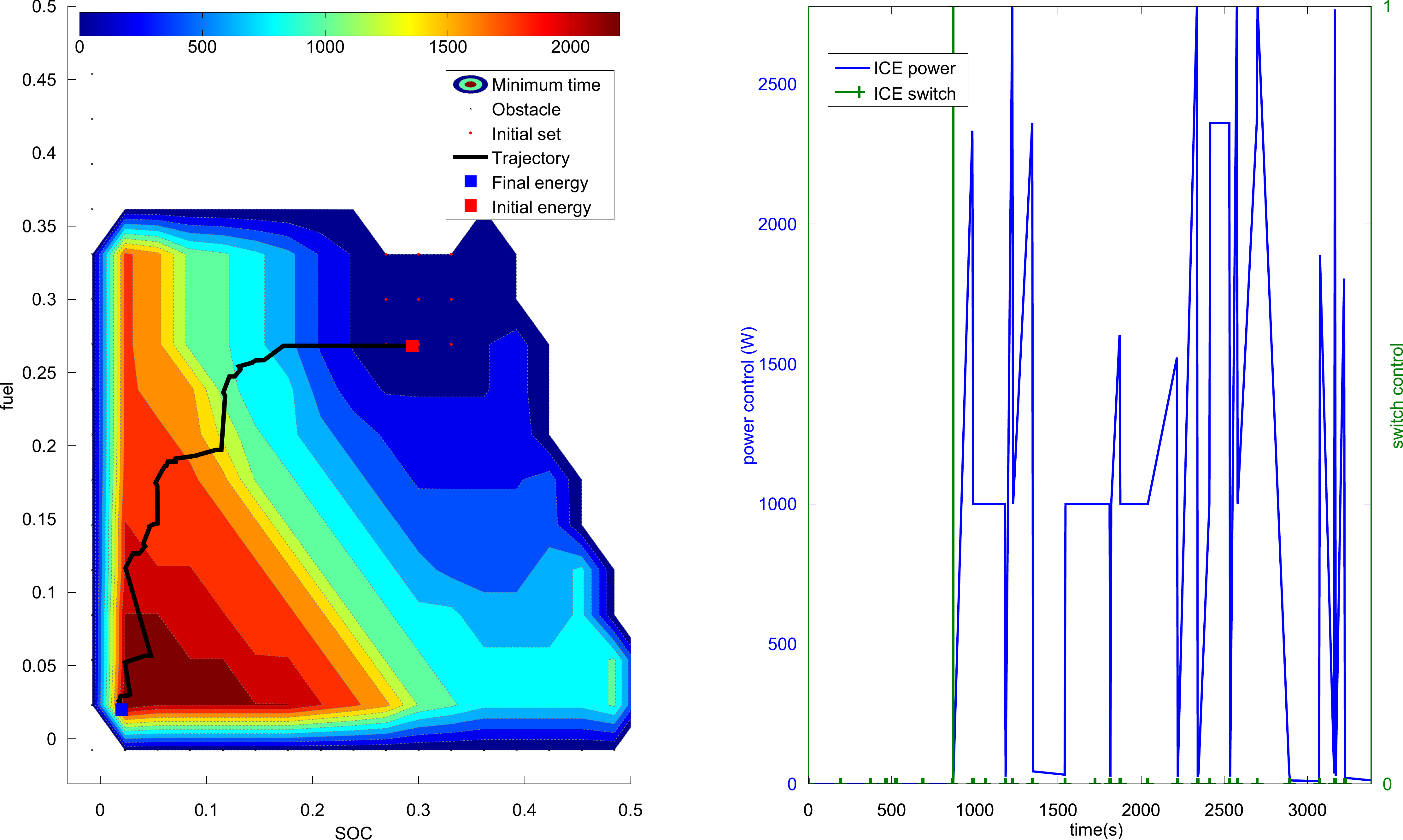}\label{vehicleAutonomyReachA} 
	\end{center}
	\caption{Figure showing the trajectory of the energy state from an valid initial state until the last reachable point.
		The plot on the left groups the contour levels of the minimum time function, the obstacle and initial sets.
		On the right, the figure displays the range extender power control $u$ and the switch sequence $w$ controlling the ICE state.
		The synthesis is obtained using step sizes of $\Delta x = 0.0308$. }
	\end{figure}

	We proceed to the computation of the autonomy using a realistic vehicle model.
	Table \ref{simulationResultsCAR} summarizes the main results.
	Figure \ref{vehicleAutonomyReach} regroups the minimum time function, the synthesized controller that allows the vehicle to reach its maximum range and the corresponding evolution of the SOC and fuel.
	The autonomies computed in both cases can be compared to the autonomy of the vehicle under a purely electric utilization.
	In the first case, the utilization of the range extender roughly doubles the range of the vehicle (a $104.7\%$ range increase).
	When considering the initial condition with a SOC at $60\%$ and fuel at $40\%$, the range extender increases the vehicle range by over $60\%$.
	We also consider the total operating \emph{financial cost} (given in euros in this study) of the range extender.
	This cost is based on an estimated $\eta = 1.5$ \euro/l of fuel.
	The range extender operating cost is then normalized with respect to $100$km to yield a more meaningful result.
	We stress that we take into consideration in the cost of the range extender operation only the \emph{additional distance} due the range extender.
	Therefore, the values of operating cost should be read as euros per 100km of range extender utilization instead of 100km of driving distance.
	In our model the reservoir has a total capacity of $6$l.

	\begin{table} \label{simulationResultsCAR}
	\centering
		\begin{tabular}{|l|c|c|} 
		\hline
		Initial SOC, fuel 			& $(0.3,0.3)$ 			& $(0.6,0.4)$			\\ \hline
		Max range 			& $45.126$ km			& $78.405$ km			\\ 
		EV range			& $22.045$ km			& $48.209$ km			\\ 	
		Relative range increase  		& $104.70 \%$			& $62.64 \%$			\\ \hline
		Range extender op. cost 	& $11.70 $ \euro$/100$km	& $11.92 $ \euro$/100$km	\\ \hline
		\end{tabular}
		\caption{Results for the autonomy of the range extender electric vehicle.
			The table presents the maximum range of the vehicle when equipped with a range extender (Max range) and the range of the same vehicle when operating in a purely electric mode (EV range).
			The financial cost of operation of the range extender is estimated, based on a price of $1.5$ \euro/l.}
	\end{table}

	\section{Conclusion}

	In this study we present an algorithm to determine the autonomy of a hybrid vehicle.
	The algorithm uses information about the driving profile of the vehicle, obtained typically through a navigation system.
	Once the destination is fixed and if there is not enough energy on-board the vehicle to reach the destination point, the algorithms computes the longest driving distance and synthesize the associated control sequence of the controllable power source.
	In our application, we focus on the range extender electric vehicle-class of hybrid vehicles.
	Our model considers an internal combustion engine as a range extender and a high voltage DC battery as main energy source.
	The dynamic evolution of the state is considered to be a discrete-time evolution.
	
	The algorithm computes the reachable set, set of reachable states, given a final time and initial state, via a level-set approach.
	The reachable set is characterized as the negative region of a value function associated to an optimal control problem.
	To compute the value function, we formulate an optimal control problem with obstacle and we state a dynamic programming principle verified by the optimal control problem.
	Then, the value function can be computed using a (deterministic) dynamic programming principle.

	Numerical simulations of this procedure are performed first using a hybrid vehicle toy model and then using a realistic range extender electric vehicle model.
	The results shows that, depending on the initial conditions, operating the range extender as described in this paper may yield a $100\%$ increase in the vehicle range, with respect to a purely electric vehicle.

\end{document}